\documentclass[12pt]{amsart}

\usepackage{amsmath}
\usepackage{amscd}
\usepackage{amssymb}
\usepackage{latexsym}
\usepackage{stmaryrd} 
\DeclareMathAlphabet{\mathbbb}{U}{bbold}{m}{n}

\newtheorem{theorem}{Theorem}[section]
\newtheorem*{theorem*}{Theorem}
\newtheorem{lemma}[theorem]{Lemma}
\newtheorem*{lemma*}{Lemma}
\newtheorem{corollary}[theorem]{Corollary}
\newtheorem{proposition}[theorem]{Proposition}

\newtheorem{remark}[theorem]{Remark}
\newtheorem{definition}[theorem]{Definition}

\makeatletter
\def\revddots{\mathinner{\mkern1mu\raise\p@
\vbox{\kern7\p@\hbox{.}}\mkern2mu
\raise4\p@\hbox{.}\mkern2mu\raise7\p@\hbox{.}\mkern1mu}}
\makeatother % ddots gespiegelt
\newcommand{\bgl}{\begin{equation}} %eine Gleichung mit Ziffer
\newcommand{\egl}{\end{equation}}
\newcommand{\bgloz}{\begin{equation*}} %eine Gleichung ohne Ziffer
\newcommand{\egloz}{\end{equation*}}
\newcommand{\bgln}{\begin{eqnarray}} %mehrere Gleichungen mit Ziffer
\newcommand{\egln}{\end{eqnarray}}
\newcommand{\bglnoz}{\begin{eqnarray*}} %mehrere Gleichungen ohne Ziffer
\newcommand{\eglnoz}{\end{eqnarray*}}
\newcommand{\btheo}{\begin{theorem}}
\newcommand{\etheo}{\end{theorem}}
\newcommand{\btheooz}{\begin{theorem*}}
\newcommand{\etheooz}{\end{theorem*}}
\newcommand{\blemma}{\begin{lemma}}
\newcommand{\elemma}{\end{lemma}}
\newcommand{\blemmaoz}{\begin{lemma*}}
\newcommand{\elemmaoz}{\end{lemma*}}
\newcommand{\bproof}{\begin{proof}}
\newcommand{\eproof}{\end{proof}}
\newcommand{\bbew}{\begin{beweis}}
\newcommand{\ebew}{\end{beweis}}
\newcommand{\bremark}{\begin{remark}\em}
\newcommand{\eremark}{\end{remark}}
\newcommand{\bdefin}{\begin{definition}}
\newcommand{\edefin}{\end{definition}}
\newcommand{\bprop}{\begin{proposition}}
\newcommand{\eprop}{\end{proposition}}
\newcommand{\bcor}{\begin{corollary}}
\newcommand{\ecor}{\end{corollary}}
\newcommand{\bfa}{\begin{cases}} %Fallunterscheidung
\newcommand{\efa}{\end{cases}}
\newcommand{\cL}{\mathcal L}

\def\Az{\mathbb{A}}
\def\Cz{\mathbb{C}}
\def\Fz{\mathbb{F}}

\def\Qz{\mathbb{Q}}

\def\Tz{\mathbb{T}}
\def\Zz{\mathbb{Z}}

\def\1z{\mathbbb{1}}
\newcommand{\fA}{\mathfrak A}

\newcommand{\an}[1]{``#1''} % Anfuehrungsstriche
\newcommand{\ti}{\tilde}

\newcommand{\ma}{\mapsto} % wird abgebildet auf
\newcommand\into{\hookrightarrow} % injektiv
\def\SEMI{\mbox{$\times\kern-2pt\vrule height5pt width.6pt \kern3pt $}}

\newcommand{\id}{{\rm id}}

\newcommand{\tei}{\mid} % teilt
\newcommand{\reg}{^\times} % regulaer
\newcommand{\nneg}{_{\geq 0}} % nicht negativ
\newcommand{\abs}[1]{\lvert#1\rvert} % Betrag
\newcommand{\defeq}{\mathrel{:=}} % per Definition
\newcommand{\dop}{\text{: }} % in Mengen
\newcommand{\fa}{\text{ for all }} % fuer alle
\newcommand{\ilim}{\varinjlim} % induktiver Limes
\newcommand{\plim}{\varprojlim} % projektiver Limes
\newcommand{\extalg}{\Lambda^* \,} % aeußeres Produkt
\newcommand{\lge}{\left\{} % links geschweift
\newcommand{\rge}{\right\}} % rechts geschweift
\newcommand{\lru}{\left(} % links rund
\newcommand{\rru}{\right)} % rechts rund
\newcommand{\leck}{\left[} % links eckig
\newcommand{\reck}{\right]} % rechts eckig
\newcommand{\lsp}{\left\langle} % links spitz
\newcommand{\rsp}{\right\rangle} % links spitz
\newcommand{\rukl}[1]{\lru #1 \rru} % runde Klammer
\newcommand{\eckl}[1]{\leck #1 \reck} % eckige Klammer
\newcommand{\gekl}[1]{\lge #1 \rge} % geschweifte Klammer
\newcommand{\spkl}[1]{\lsp #1 \rsp} % spitze Klammer
\newcommand{\menge}[2]{\gekl{ #1 \dop #2 }} % Menge
\newcommand{\sig}{^{(\sigma)}} % regulaer
\newcommand{\dualF}{\widehat{\Fz_q \reg}} % duale Gruppe der multiplikativen Gruppe von F_q 
\begin{document}

\title[K-theory for ring C*-algebras attached to function fields]{K-theory for ring C*-algebras attached to function fields with only one infinite place}

\author{Xin Li}


\thanks{\scriptsize{Research supported by the Deutsche Forschungsgemeinschaft (SFB 878).}}

\begin{abstract}
We study the K-theory of ring C*-algebras associated to rings of integers in global function fields with only one single infinite place. First, we compute the torsion-free part of the K-groups of these ring C*-algebras. Secondly, we show that, under a certain primeness condition, the torsion part of K-theory determines the inertia degrees at infinity of our function fields.
\end{abstract}

\maketitle

\tableofcontents

\setlength{\parindent}{0pt} \setlength{\parskip}{0.5cm}

\section{Introduction}

We generalize the work of J. Cuntz and the author in \cite{Cu-Li3}, where K-theory has been computed for ring C*-algebras attached to $\Fz_q[T]$, the polynomial ring over the finite field with $q$ elements. Here $q$ is a prime power, i.e. $q=p^\nu$ for some prime $p$ and positive integer $\nu$. The construction of ring C*-algebras is recalled in the next section.

In the present paper, we consider more general global function fields, i.e. finite separable field extensions of $\Fz_p(T)$, namely those which have only one single infinite place. In other words, we focus on function fields that do not split at infinity. Here, our convention is to call the place of $\Fz_p(T)$ corresponding to the inverse of the indeterminant $T$ the infinite place (as it is usually done). A place of a function field is then called infinite if it sits above this infinite place of $\Fz_p(T)$.

Given any function field, we can consider the ring C*-algebra attached to its ring of integers. By the ring of integers, we mean the integral closure of $\Fz_p[T]$ in this function field. Our present goal is to study the K-theory of this ring C*-algebra in the situation of a function field with only one single infinite place.

In our previous work \cite{Cu-Li3}, we have considered specific examples, namely the ring of integers $\Fz_q[T]$ in the function field $\Fz_q(T) / \Fz_p(T)$. So the setting in \cite{Cu-Li3} was very concrete. On the one hand, this had the advantage that the K-theoretic computation could be done in a very explicit way. This led to explicit generators for the K-groups of interest. On the other hand, it was not clear how to generalize the work in \cite{Cu-Li3}, i.e. what happens in the case of more general function fields.

The present piece of work clarifies the situation and helps to understand what is going on in the general case. Our investigations reveal which properties of the function field are needed for the approach developed in \cite{Cu-Li3} to work out, and where the obstacles lie to more general K-theoretic computations.

Moreover, new phenomena arise in the more general situation. Most notably, we will see that the K-theory of ring C*-algebras attached to function fields (or rather their rings of integers) can carry torsion. This is not the case for the ring C*-algebra of $\Fz_q[T]$ as shown in \cite{Cu-Li3}. Unfortunately, the appearance of torsion leads to extension problems in our computations that cannot be solved in general. For this reason, we can in general only determine the torsion-free part of the K-groups. However, it is at the same time an interesting problem to find out what sort of information about the function fields can be extracted from the torsion part of the K-groups of the associated ring C*-algebras. It turns out that torsion in K-theory encodes the inertia degrees at infinity of our function fields.

Let us now state our main results:

Let $K / \Fz_p(T)$ be a finite separable field extension. Assume that $\Fz_q$ is the field of constants in $K$, i.e. $\Fz_q$ is the algebraic closure of $\Fz_p$ in $K$. We assume that $K$ has only one single infinite place. Let $f$ be the inertia degree at infinity of $K / \Fz_q(T)$. Define $Q^f \defeq \tfrac{q^f-1}{q-1}$. Moreover, the multiplicative group $K\reg$ always admits a decomposition $K\reg = \Fz_q\reg \times \Gamma$ where $\Gamma$ is a free abelian subgroup of $K\reg$. Now let $R$ be the ring of integers of $K$. In this situation, we have

\btheo
\label{mainthm1}
After inverting $Q^f$, we obtain for the K-theory of $\fA[R]$, the ring C*-algebra of $R$:
\bgloz
  K_*(\fA[R]) \otimes_{\Zz} \Zz[\tfrac{1}{Q^f}] \cong \rukl{\ti{K}_0(C^*(\Fz_q\reg)) \otimes_{\Zz} \extalg(\Gamma)} \otimes_{\Zz} \Zz[\tfrac{1}{Q^f}].
\egloz
\etheo
$\ti{K}_0$ stands for reduced K-theory, and we know that $\ti{K}_0(C^*(\Fz_q\reg)) \cong \Zz^{q-2}$. $\extalg(\Gamma)$ is the exterior $\Zz$-algebra over $\Gamma$. Moreover, $K_*$ denotes the direct sum of $K_0$ and $K_1$ and is viewed as a $\Zz / 2 \Zz$-graded abelian group. So we are identifying $\Zz / 2 \Zz$-graded abelian groups. On $\Zz[\tfrac{1}{Q^f}]$ and $\ti{K}_0(C^*(\Fz_q\reg))$, we take the trivial grading, and on $\extalg(\Gamma)$ we take the canonical grading. All the tensor products are graded tensor products.

For our second main result, take two function fields $K_1$ and $K_2$. We assume that both of them have the field of constants $\Fz_q$ and only one single infinite place. Let $f_i$ be the inertia degree at infinity of $K_i / \Fz_q(T)$ for $i=1,2$. $R_i$ again denotes the ring of integers in $K_i / \Fz_q(T)$ for $i=1,2$. Under the additional assumption that $(f_i,q-1)=1$ for $i=1,2$, we obtain
\btheo
\label{mainthm2}
If $K_*(\fA[R_1]) \cong K_*(\fA[R_2])$, then the inertia degrees at infinity of the extensions $K_1 / \Fz_q(T)$ and $K_2 / \Fz_q(T)$ coincide, i.e. we have $f_1=f_2$.
\etheo

This paper is organized as follows: In Section~\ref{pre}, we explain the construction of ring C*-algebras for rings of integers in function fields. Moreover, we review the work in \cite{Cu-Li3} in order to recall the strategy of the K-theoretic computations in the case of polynomial rings over finite fields. As in \cite{Cu-Li3}, it is more convenient to work over the finite place of $\Fz_p(T)$ corresponding to the indeterminant $T$ itself instead of the infinite place corresponding to $T^{-1}$. We explain the necessary reduction step in Section~\ref{inf-fin}. Our original task can be reformulated as computing K-theory for a certain crossed product by the $ax+b$-group over a function field. We start with the action of the additive group and of the roots of unity in Section~\ref{add&rou}. It then remains to treat the full multiplicative group. This is achieved in two steps: The first step is to adjoin a prime element (Section~\ref{prime}). The second step is to add the remaining part of the multiplicative group. We show that it is at least possible to do this step rationally (Section~\ref{rational}), so that we obtain a proof of Theorem~\ref{mainthm1}. Moreover, under a certain primeness condition, it is possible to say a bit more about the torsion part of K-theory. This observation, together with a little number theoretic lemma, gives a proof of Theorem~\ref{mainthm2} (Section~\ref{torsion}). 

I would like to thank U. Hartl for helpful comments about global function fields.

\section{Preliminaries}
\label{pre}

\subsection{Ring C*-algebras}
The theory of ring C*-algebras has been initiated in \cite{Cun}. From that point of departure, it has been further developed in \cite{Cu-Li1}, \cite{Li} and \cite{Cu-Li2}. We refer the reader to the references mentioned above for the general theory of ring C*-algebras. Now, let us explain the construction in the concrete case of function fields.

By a (global) function field, we mean a finite separable field extension $K / \Fz_p(T)$. Here $p$ is a prime number and $\Fz_p(T)$ is the quotient field of the polynomial ring over the finite field with $p$ elements. The ring of integers $R$ in $K / \Fz_p(T)$ is the integral closure of $\Fz_p[T]$ in $K$, i.e. $R$ is the following subring of $K$:
\bgloz
  \menge{x \in K}{\exists \: a_0, \dotsc, a_{m-1} \in \Fz_p[T] \text{ with } x^m + a_{m-1}x^{m-1} + \dotsb + a_1x + a_0 = 0}.
\egloz
The reader may consult \cite{Ro} for an introduction to the theory of function fields. The theory of places and adeles can be found in \cite{Weil}, Chapters~I-IV.

To construct the ring C*-algebra of $R$, consider the Hilbert space $\ell^2(R)$ of square integrable functions on $R$. Denote by $\menge{\varepsilon_r}{r \in R}$ the canonical orthonormal basis of $\ell^2(R)$, i.e. $\varepsilon_r(s) = \delta_{r,s}$ for $s$ in $R$ ($\delta_{r,s} = 1$ if $r=s$ and $\delta_{r,s} = 0$ if $r \neq s$). The ring structure of $R$ allows us to construct two natural families of bounded operators on $\ell^2(R)$, namely
\begin{itemize}
  \item
  additive shifts determined by $U^b \varepsilon_r = \varepsilon_{b+r}$ for every $b \in R$, $r \in R$ 
  \item
  multiplicative shifts given by $S_a \varepsilon_r = \varepsilon_{ar}$ for all $0 \neq a \in R$, $r \in R$.
\end{itemize}
The ring C*-algebra of $R$ is then given by the C*-subalgebra of $\cL(\ell^2(R))$ generated by $\menge{U^b}{b \in R}$ and $\menge{S_a}{0 \neq a \in R}$. It is denoted by $\fA[R]$, i.e. we have
\bgloz
  \fA[R] \defeq C^*(\menge{U^b}{b \in R} \cup \menge{S_a}{0 \neq a \in R}) \subseteq \cL(\ell^2(R)).
\egloz
In our situation of function fields, the C*-algebra $\fA[R]$ is purely infinite and simple (compare \cite{Cu-Li1}, Theorem~3.6). Moreover, it admits -- up to Morita equivalence -- a crossed product description
\bgl
\label{A-Mor-finadele}
  \fA[R] \sim_M C_0(\Az_f) \rtimes K \rtimes K\reg.
\egl
This follows more or less directly from our constructions (see \cite{Cu-Li1}, Section~5). Here $\Az_f$ is the finite adele ring over $K$. Our function field $K$ embeds diagonally into $\Az_f$, so there is a canonical action of the $ax+b$-group $K \rtimes K\reg$ over $K$ on $C_0(\Az_f)$ by affine transformations. This is precisely the action which gives rise to the crossed product on the right hand side of \eqref{A-Mor-finadele}.

\subsection{Previous K-theoretic computations}
The K-theory of $\fA[R]$ has been computed in the case of the function field $\Fz_q(T) / \Fz_p(T)$ in \cite{Cu-Li3}, where $q$ is a power of the prime number $p$ and we take the canonical embedding $\Fz_p(T) \into \Fz_q(T)$. For this function field $\Fz_q(T) / \Fz_p(T)$, the ring of integers $R$ is the polynomial ring with coefficients in $\Fz_q$, i.e. $R = \Fz_q[T]$. The main result in \cite{Cu-Li3} says that the K-theory of $\fA[\Fz_q[T]]$ can be described as follows:

The multiplicative group $\Fz_q(T)\reg$ admits a decomposition $\Fz_q(T)\reg = \Fz_q\reg \times \Gamma$ where $\Gamma$ is a free abelian subgroup of $\Fz_q(T)\reg$ ($\Gamma$ is actually a free abelian group on countably infinitely many generators). Theorem~7.2 in \cite{Cu-Li3} tells us that
\bgl
\label{K(A[F_q[T]])}
  K_*(\fA [\Fz_q[T]]) \cong \ti{K}_0(C^*(\Fz_q\reg)) \otimes_{\Zz} \extalg(\Gamma)
\egl
as $\Zz / 2 \Zz$-graded abelian groups. Here we use the same notations as in Theorem~\ref{mainthm1}.

In the sequel, let us briefly recall the steps that lead to \eqref{K(A[F_q[T]])}.

The first and most important ingredient is the so-called duality theorem, proven in \cite{Cu-Li2}. It says the following:

Let $K$ be a global field, i.e. we have a finite seperable field extension $K / \Qz$ or $K / \Fz_p(T)$. We can associate to $K$ its finite and its infinite adele space, denoted by $\Az_f$ and $\Az_{\infty}$, respectively. Since $K$ embeds diagonally into $\Az_f$ and $\Az_{\infty}$, there are canonical actions of the $ax+b$-group $K \rtimes K\reg$ on $C_0(\Az_f)$ and $C_0(\Az_{\infty})$ by affine transformations. Corollary~4.2 in \cite{Cu-Li2} tells us that the corresponding crossed products are Morita equivalent, i.e.
\bgl
\label{duality-thm}
  C_0(\Az_f) \rtimes K \rtimes K\reg \sim_M C_0(\Az_{\infty}) \rtimes K \rtimes K\reg.
\egl
Of course, a function field is a special case of a global field, so we can apply \eqref{duality-thm} to function fields. Combining this with \eqref{A-Mor-finadele}, we obtain for every function field $K / \Fz_p(T)$ with ring of integers $R$ that
\bgl
\label{A-Mor-infadele}
  \fA[R] \sim_M C_0(\Az_{\infty}) \rtimes K \rtimes K\reg.
\egl
This is the first step in our computation; it reduces the K-theory computation of $\fA[R]$ to the computation of the K-theory of $C_0(\Az_{\infty}) \rtimes K \rtimes K\reg$.

The second step, another reduction step, is to switch places from the infinite place to the finite place corresponding to the indeterminant $T$ in $\Fz_q(T)$. For the function field $\Fz_q(T) / \Fz_p(T)$ treated in \cite{Cu-Li3}, this means the following: The infinite adele space $\Az_{\infty}$ of $\Fz_q(T) / \Fz_p(T)$ is given by the ring of formal Laurent series $\Fz_q((T))$ with the embedding $\Fz_q(T) \ni f(T) \ma f(T^{-1}) \in \Fz_q((T))$, i.e. the indeterminant is replaced by its inverse. The local field of the place of $\Fz_q(T)$ corresponding to the indeterminant $T$ is again given by $\Fz_q((T))$, but this time with the canonical embedding $\Fz_q(T) \ni f(T) \ma f(T) \in \Fz_q((T))$. These two embeddings give rise to two actions of $\Fz_q(T) \rtimes \Fz_q(T)\reg$ on $C_0(\Az_{\infty})$ and on $C_0(\Fz_q((T)))$. Thus we obtain two crossed products. In Section~3 of \cite{Cu-Li3}, we prove that
\bgloz
  C_0(\Az_{\infty}) \rtimes \Fz_q(T) \rtimes \Fz_q(T)\reg \cong C_0(\Fz_q((T))) \rtimes \Fz_q(T) \rtimes \Fz_q(T)\reg.
\egloz
Thus to compute the K-theory of $\fA[\Fz_q[T]]$, we have to determine the K-theory of $C_0(\Fz_q((T))) \rtimes \Fz_q(T) \rtimes \Fz_q(T)\reg$. This computation is then done in three steps:

As a first step, we treat the additive group and the roots of unity, i.e. we compute K-theory for $C_0(\Fz_q((T))) \rtimes \Fz_q(T) \rtimes \Fz_q\reg$ (see \cite{Cu-Li3}, Corollary~5.9; $C_0(\Fz_q((T))) \rtimes \Fz_q(T) \rtimes \Fz_q\reg$ is abbreviated by $A_{-1}$ in \cite{Cu-Li3}). The idea is to find suitable filtrations and to use continuity of K-theory.

The second step is to adjoin the canonical prime element for the place corresponding to $T$, namely the indeterminant $T$ itself. This means that we have to compute K-theory for $C_0(\Fz_q((T))) \rtimes \Fz_q(T) \rtimes (\Fz_q\reg \times \spkl{T})$. However, it is crucial for what comes next that we do not only compute K-theory, but even obtain explicit generators with a certain invariance property (see \cite{Cu-Li3}, Proposition~5.10 and Lemma~6.1; $C_0(\Fz_q((T))) \rtimes \Fz_q(T) \rtimes (\Fz_q\reg \times \spkl{T})$ is abbreviated by $A_0$ in \cite{Cu-Li3}). This explicit K-theoretic computation is achieved through a detailed analysis of the Pimsner-Voiculescu sequence applied to our situation.

Finally, the third step is to adjoin the remaining part of the multiplicative group $\Fz_q(T)\reg$. Here, it turns out to be crucial that we can find generators for this remaining part of $\Fz_q(T)\reg$ which all lie in the coset $1+T\Fz_q[T]$, i.e. which are congruent to $1$ modulo $T\Fz_q[T]$. The idea in this step is to use the Pimsner-Voiculescu sequence and to compare our situation to the one of higher dimensional tori. This comparison works out precisely because of the invariance property of the generators for K-theory from the previous step.

This is how the K-theory of $\fA[\Fz_q[T]]$ can be computed. In the present paper, we follow the same route. Some steps work out just as in the case of $\Fz_q(T) / \Fz_p(T)$, but sometimes we need modifications or new ideas to make things work.

\bremark
K-theory for ring C*-algebras has also been computed in the case of number fields, see \cite{Cu-Li2}. As for function fields, the duality theorem \eqref{duality-thm} plays an important role in the computations for number fields, but for quite different reasons.
\eremark

\section{Switching places from infinite to finite}
\label{inf-fin}

Let now $K / \Fz_p(T)$ be a function field with field of constants $\Fz_q$. This means that $\Fz_q$ is the algebraic closure of $\Fz_p$ in $K$. Let $R$ be the ring of integers in $K / \Fz_p(T)$. We know by \eqref{A-Mor-infadele} that the ring C*-algebra of $R$ is Morita equivalent to $C_0(\Az_{\infty}) \rtimes K \rtimes K\reg$. As in the case of polynomial rings and their quotient fields, it is more convenient to work not with the infinite place $v_{\infty}$ but with the place $v_T$ of $\Fz_p(T)$ corresponding to the prime ideal of $\Fz_p[T]$ generated by the indeterminant $T$. This means that we would like to work with the places of $K / \Fz_p(T)$ sitting above $v_T$ instead of those places which sit above $v_{\infty}$. It is possible to arrange this by twisting the field extension $K / \Fz_p(T)$ by an automorphism of the ground field.

Let $\iota$ be the inclusion $\Fz_p(T) \into K$ corresponding to the given field extension $K / \Fz_p(T)$. Consider the automorphism $\sigma$ of $\Fz_p(T)$ which fixes $\Fz_p$ and sends $T$ to $T^{-1}$. Compose $\iota$ with $\sigma$ to obtain a new field extension $K\sig / \Fz_p(T)$. In other words, we take the same field $K$, but view $\Fz_p(T)$ as a subfield via $\iota \circ \sigma$ and not via $\iota$ as it was the case for the extension $K / \Fz_p(T)$. Since $K\sig = K$ as fields, $K\sig$ has the same field of constants (namely $\Fz_q$) and the same places as $K$. Moreover, it is clear from our construction that a place $w$ of $K\sig$ lies above the place $v_T$ of $\Fz_p(T) \overset{\iota \circ \sigma}{\into} K\sig$ if and only if the same place $w$, now viewed as a place of $K$, sits above the infinite place $v_{\infty}$ of $\Fz_p(T) \overset{\iota}{\into} K$. Therefore we have
\bgl
\label{A_T=A_inf}
  \prod_{w \tei v_T} K_w\sig = \prod_{w \tei v_{\infty}} K_w.
\egl
Let $\Az_T\sig = \prod_{w \tei v_T} K_w\sig$. The right hand side of \eqref{A_T=A_inf} is $\Az_{\infty}$, the infinite adele space of $K$, so that we obtain
\bgl
\label{AKK-Mor-AKsigKsig}
  C_0(\Az_{\infty}) \rtimes K \rtimes K\reg = C_0(\Az_T\sig) \rtimes K\sig \rtimes (K\sig)\reg.
\egl
Therefore, our goal is to compute the K-theory of $C_0(\Az_T\sig) \rtimes K\sig \rtimes (K\sig)\reg$.

In analogy to the case of $K / \Fz_p(T)$, we denote by $R\sig$ the ring of integers of the function field $K\sig / \Fz_p(T)$. Let $\overline{R\sig}$ be the closure of $R\sig$ in $\Az_T\sig$, i.e. 
\bgloz
  \overline{R\sig} = \prod_{w \tei v_T} R_w\sig,
\egloz
where $R_w\sig$ is the closure of $R\sig$ in $K_w\sig$.

\section{Additive group and roots of unity}
\label{add&rou}
Let us first of all state the main result of this section. For non-negative integers $m$, let $\1z_m$ denote the characteristic function of $T^m \overline{R\sig} \subseteq \Az_T\sig$. As $T^m \overline{R\sig}$ is closed and open in $\Az_T\sig$, $\1z_m$ lies in $C_0(\Az_T\sig)$. In the sequel, we will often deal with crossed products by the $ax+b$-group $K\sig \rtimes (K\sig)\reg$ or by subgroups of this group. In each of these crossed products, let $v^b$ denote the unitaries in the multiplier algebra of the crossed product corresponding to the additive part, and let $t_a$ denote the unitaries for the multiplicative part. Moreover, whenever in such a crossed product by (a subgroup of) $K\sig \rtimes (K\sig)\reg$, $\Fz_q\reg$ is contained in the multiplicative part of the group, then for a character $\chi \in \dualF$ let $p_\chi$ be the spectral projection $\tfrac{1}{q-1} \sum_{a \in \Fz_q\reg} \chi(a) t_a$ in the multiplier algebra of the crossed product.

Here is the main result of this section:

Let $n$ be the degree of $K\sig$ over $\Fz_q(T)$. Moreover, $1 \in \dualF$ denotes the trivial character.
\bprop
\label{K(add&rou)}
We can identify $K_0(C_0(\Az_T\sig) \rtimes K\sig \rtimes \Fz_q\reg)$ with 
\bgloz
  \Zz[\tfrac{1}{q}] \oplus \bigoplus_{\dualF \setminus \gekl{1}} \Zz
\egloz
in such a way that the $K_0$-class $[\1z_m]$ corresponds to 
$
  \rukl{
  \begin{smallmatrix}
  q^{-mn} \\
  0 \\
  \vdots \\
  0
  \end{smallmatrix}
  } 
$
and the $K_0$-class $[\1z_m p_\chi]$ corresponds to
$
  \rukl{
  \begin{smallmatrix}
  - q^{-mn} Q^{mn}\\
  0 \\
  \vdots \\
  0 \\
  1 \\
  0 \\
  \vdots \\
  0
  \end{smallmatrix}
  }
$
for every $m \in \Zz \nneg$ and $1 \neq \chi \in \dualF$. Here the entry \an{1} in the vector corresponding to $[\1z_m p_\chi]$ is the $\chi$-th entry in $\bigoplus_{\dualF \setminus \gekl{1}} \Zz$. Moreover, we write $Q^{mn}$ for $\tfrac{q^{mn}-1}{q-1}$.

$K_1(C_0(\Az_T\sig) \rtimes K\sig \rtimes \Fz_q\reg)$ vanishes.
\eprop

The first step in the proof of this proposition is to compute the K-theory of $C(\overline{R\sig}) \rtimes R\sig \rtimes \Fz_q\reg$. Since $K / \Fz_q(T)$ is separable, $K\sig / \Fz_q(T)$ is separable as well. In addition, $\Fz_q[T]$ is a principal ideal domain. Therefore, there is $\Fz_q[T]$-basis for $R\sig$ (see \cite{Neu}, Chapter~I, Proposition~(2.10)). Let $\omega_1, \dotsc, \omega_n$ be such a basis. We set
\bgloz
  R_m\sig \defeq \menge{\sum_{i=1}^n f_i \omega_i}{f_i \in \Fz_q[T], \deg(f_i) \leq m-1}.
\egloz
We know that $\overline{R\sig}$ can be identified with $\plim_m \gekl{R\sig / T^m R\sig; p_{m+1,m}}$ where $p_{m+1,m}: R\sig / T^{m+1} R\sig \to R\sig / T^m R\sig$ is the canonical projection. Moreover, we have $R\sig = \bigcup_{m=1}^{\infty} R_m\sig$. It follows that we can identify $C(\overline{R\sig}) \rtimes R\sig \rtimes \Fz_q\reg$ with
\bgloz
  \ilim_m \gekl{C(R\sig / T^m R\sig) \rtimes R_m\sig \rtimes \Fz_q\reg; i_{m,m+1}}
\egloz
where $i_{m,m+1}: C(R\sig / T^m R\sig) \rtimes R_m\sig \rtimes \Fz_q\reg \to C(R\sig / T^{m+1} R\sig) \rtimes R_{m+1}\sig \rtimes \Fz_q\reg$ is given by
\bgl
\label{i_m,m+1}
  i_{m,m+1}(g v^b t_a) = (g \circ p_{m+1,m}) v^b t_a
\egl
for all $g \in C(R\sig / T^m R\sig)$, $b \in R_m\sig$ and $a \in \Fz_q\reg$. This can be proven analogously to Lemma~5.3 in \cite{Cu-Li3}.

As a first step, we describe $C(R\sig / T^m R\sig) \rtimes R_m\sig \rtimes \Fz_q\reg$. Let $V_a$, for $a \in \Fz_q \reg$, be the canonical unitary generators of $C^*(\Fz_q \reg)$. Moreover, let $e_m \in C(R\sig / T^m R\sig)$ be the function with value $1$ at $0 + T^m R\sig$ and $0$ everywhere else. And for every $d$, $d'$ in $R_m\sig$, let $e_{d,d'}$ be the rank one operator in $\cL(\ell^2(R\sig / T^m R\sig))$ given by $e_{d,d'}(\xi) = \spkl{\xi, \varepsilon_{d'+T^m R\sig}} \varepsilon_{d+T^m R\sig}$. Here $\gekl{\varepsilon_{b+T^m R\sig}}$ is the canonical orthonormal basis of $\ell^2(R\sig / T^m R\sig)$.
\blemma
We have
\bgloz
  C(R\sig / T^m R\sig) \rtimes R_m\sig \rtimes \Fz_q\reg \cong \cL(\ell^2(R\sig / T^m R\sig)) \otimes C^*(\Fz_q\reg)
\egloz
via
\bgloz
  v^b e_m v^{-b'} t_a \ma e_{b,a^{-1}b'} \otimes V_a.
\egloz
\elemma
The proof of this lemma is analogous to the one of Lemma~5.4 in \cite{Cu-Li3}.

Identifying $\cL(\ell^2(R\sig / T^m R\sig))$ with $M_{q^{mn}}(\Cz)$, this lemma leads to
\bcor
\bgloz
  K_0(C(R\sig / T^m R\sig) \rtimes R_m\sig \rtimes \Fz_q\reg) \cong \bigoplus_{\dualF} \Zz
\egloz
and free generators for $K_0$ are given by the $K_0$-classes $\menge{[e_m p_\chi]}{\chi \in \dualF}$.

Moreover, $K_1(C(R\sig / T^m R\sig) \rtimes R_m\sig \rtimes \Fz_q\reg)$ vanishes.
\ecor
This already shows, using continuity of $K_1$, that $K_1(C(\overline{R\sig}) \rtimes R\sig \rtimes \Fz_q\reg)$ vanishes.

To determine $K_0(C(\overline{R\sig}) \rtimes R\sig \rtimes \Fz_q\reg)$, we have to compute what $i_{m,m+1}$ (given by \eqref{i_m,m+1}) does on $K_0$.
\blemma
\label{i_*}
With respect to the generators
\bgloz
  \menge{[e_m p_\chi]}{\chi \in \dualF} \text{ and } \menge{[e_{m+1} p_\chi]}{\chi \in \dualF}
\egloz
of $K_0(C(R\sig / T^m R\sig) \rtimes R_m\sig \rtimes \Fz_q\reg)$ and $K_0(C(R\sig / T^{m+1} R\sig) \rtimes R_{m+1}\sig \rtimes \Fz_q\reg)$, respectively, $(i_{m,m+1})_*$ is of the form
\bgloz
  (i_{m,m+1})_* 
  = 
  \rukl{
  \begin{smallmatrix}
  1 & & 0 \\
   & \ddots & \\
  0 & & 1
  \end{smallmatrix}
  }
  +
  \rukl{
  \begin{smallmatrix}
  Q^n & \dotso & Q^n \\
  \vdots & Q^n & \vdots \\
  Q^n & \dotso & Q^n
  \end{smallmatrix}
  }. 
\egloz
As above, $Q^n = \tfrac{q^n-1}{q-1}$.
\elemma
\bproof
We have to determine the $K_0$-class of
\bgloz
  i_{m,m+1}(e_m p_\chi) = e_m p_\chi = \rukl{\sum_{b \in R\sig / T R\sig} v^{T^mb} e_{m+1} v^{-T^mb}} p_\chi
\egloz
in $K_0(C(R\sig / T^{m+1} R\sig) \rtimes R_{m+1}\sig \rtimes \Fz_q\reg)$. To do so, we first construct for every $\psi$, $\chi$ in $\dualF$ a partial isometry in $C(R\sig / T^{m+1} R\sig) \rtimes R_{m+1}\sig \rtimes \Fz_q\reg$ with support projection $(e_m-e_{m+1})p_{\overline{\psi}\chi}$ and range projection $(e_m-e_{m+1})p_\chi$. Let $\psi$ be a character in $\dualF$. Moreover, set $(R\sig / T R\sig)\reg \defeq (R\sig / T R\sig) \setminus \gekl{0 + T R\sig}$. The partial isometry mentioned above is constructed with the help of a map
\bgloz
  \Psi: (R\sig / T R\sig)\reg \to \Tz=\menge{z \in \Cz}{\abs{z}=1}
\egloz
satisfying the relation
\bgl
\label{Psi-psi}
  \Psi(a \cdot b + T R\sig) = \psi(a) \Psi(b + T R\sig)
\egl
for all $a \in \Fz_q\reg$ and $b + T R\sig \in (R\sig / T R\sig)\reg$. Let us assume that such a function exists. Take such a function $\Psi$ and set
\bgloz
  x_{\Psi}^{(m+1)} \defeq \sum_{b + T R\sig \in (R\sig / T R\sig)\reg} \Psi(b + T R\sig) v^{T^mb} e_{m+1} v^{-T^mb}.
\egloz
We have
\bgln
\label{xx}
  && x_{\Psi}^{(m+1)} (x_{\Psi}^{(m+1)})^* = (x_{\Psi}^{(m+1)})^* x_{\Psi}^{(m+1)} \\
  &=& \sum_{b + T R\sig \in (R\sig / T R\sig)\reg} v^{T^mb} e_{m+1} v^{-T^mb} = e_m-e_{m+1} \nonumber
\egln
and for every $\chi \in \dualF$:
\bgln
\label{px=xp}
  && p_\chi x_{\Psi}^{(m+1)} \\
  &=& (\tfrac{1}{q-1} \sum_{a \in \Fz_q\reg} \chi(a) t_a) (\sum_{b + T R\sig \in (R\sig / T R\sig)\reg} \Psi(b + T R\sig) v^{T^mb} e_{m+1} v^{-T^mb}) \nonumber \\
  &=& \tfrac{1}{q-1} \sum_{a,b} \chi(a) \Psi(b + T R\sig) v^{T^mab} e_{m+1} v^{-T^mab} t_a \nonumber \\
  &=& \tfrac{1}{q-1} \sum_{a,b} \psi(a) \Psi(b + T R\sig) v^{T^mab} e_{m+1} v^{-T^mab} \overline{\psi}(a) \chi(a) t_a \nonumber \\
  &\overset{\eqref{Psi-psi}}{=}& \sum_{a,b} \Psi(ab + T R\sig) v^{T^mab} e_{m+1} v^{-T^mab} \tfrac{1}{q-1} \overline{\psi}(a) \chi(a) t_a \nonumber \\
  &=& \sum_{a \in \Fz_q\reg} x_{\Psi}^{(m+1)} \tfrac{1}{q-1} \overline{\psi}(a) \chi(a) t_a \nonumber \\
  &=& x_{\Psi}^{(m+1)} p_{\overline{\psi}\chi}. \nonumber
\egln
The partial isometry with the desired properties is then given by $p_\chi x_{\Psi}^{(m+1)}$ since
\bglnoz
  && (p_\chi x_{\Psi}^{(m+1)})^* (p_\chi x_{\Psi}^{(m+1)}) 
  \overset{\eqref{px=xp}}{=} p_{\overline{\psi}\chi} (x_{\Psi}^{(m+1)})^* x_{\Psi}^{(m+1)} p_{\overline{\psi}\chi} \\
  &\overset{\eqref{xx}}{=}& p_{\overline{\psi}\chi} (e_m-e_{m+1}) p_{\overline{\psi}\chi} = (e_m-e_{m+1})p_{\overline{\psi}\chi}
\eglnoz
and
\bglnoz
  && (p_\chi x_{\Psi}^{(m+1)}) (p_\chi x_{\Psi}^{(m+1)})^* = p_\chi x_{\Psi}^{(m+1)} (x_{\Psi}^{(m+1)})^* p_\chi \\
  &\overset{\eqref{xx}}{=}& p_\chi (e_m-e_{m+1})p_\chi = (e_m-e_{m+1})p_\chi.
\eglnoz
This implies that the $K_0$-classes of $(e_m-e_{m+1})p_{\overline{\psi}\chi}$ and $(e_m-e_{m+1})p_\chi$ coincide, i.e. we have for every $\psi$, $\chi$ in $\dualF$
\bgl
\label{ep}
  \eckl{(e_m-e_{m+1})p_{\overline{\psi}\chi}} = \eckl{(e_m-e_{m+1})p_\chi} \text{ in } K_0(C(R\sig / T^{m+1} R\sig) \rtimes R_{m+1}\sig \rtimes \Fz_q\reg)).
\egl
With \eqref{ep} in mind, we can now compute
\bgl
\label{i(e)_1}
  \eckl{i_{m,m+1}(e_mp_\chi)} = \eckl{e_mp_\chi} = \eckl{e_{m+1}p_\chi} + \eckl{(e_m-e_{m+1})p_\chi}
\egl
and
\bgln
\label{i(e)_2}
  && (q-1)\eckl{(e_m-e_{m+1})p_\chi} \\
  &\overset{\eqref{ep}}{=}& \sum_{\psi \in \dualF} \eckl{(e_m-e_{m+1})p_{\overline{\psi}\chi}} = \eckl{e_m-e_{m+1}} \\
  &=& (q^n-1)\eckl{e_{m+1}} = (q^n-1)\sum_{\psi \in \dualF} \eckl{e_{m+1}p_\psi} \nonumber.
\egln
\eqref{i(e)_2} implies that for every $\chi$ in $\dualF$, we have
\bgl
\label{i(e)_3}
  \eckl{(e_m-e_{m+1})p_\chi} = Q^n \sum_{\psi \in \dualF} \eckl{e_{m+1} p_\psi}.
\egl
Inserting \eqref{i(e)_3} in \eqref{i(e)_1}, we obtain
\bgloz
  \eckl{i_{m,m+1}(e_mp_\chi)} = \eckl{e_{m+1}p_\chi} + Q^n \sum_{\psi \in \dualF} \eckl{e_{m+1} p_\psi}.
\egloz
Therefore, $i_{m,m+1}$ is indeed of the desired form.

However, to make all this work, we still need to construct a function $\Psi$ with property \eqref{Psi-psi}.
\blemma
\label{Psi-exists}
For every character $\psi$ in $\dualF$ there exists a function
\bgloz
  \Psi: (R\sig / T R\sig)\reg \to \Tz=\menge{z \in \Cz}{\abs{z}=1}
\egloz
satisfying the relation \eqref{Psi-psi}, i.e.
\bgloz
  \Psi(a \cdot b + T R\sig) = \psi(a) \Psi(b + T R\sig) \fa a \in \Fz_q\reg, b + T R\sig \in (R\sig / T R\sig)\reg.
\egloz
\elemma
\bproof
Let $\omega_1, \dotsc, \omega_n$ be a $\Fz_q[T]$-basis of $R\sig$ (see the beginning of this section, right after Proposition~\ref{K(add&rou)}). Then a complete system of representatives in $R\sig / T R\sig$ is given by $R_1\sig = \menge{\sum_{i=1}^n b_i \omega_i}{b_i \in \Fz_q}$. Of course, we have $\sum_{i=1}^n b_i \omega_i + T R\sig = 0 + T R\sig$ if and only if $b_1 = \dotsb = b_n = 0$. Thus every $b + T R\sig$ in $(R\sig / T R\sig)\reg$ is of the form $\sum_{i=i_0}^n b_i \omega_i$ with an index $i_0 \in \gekl{1, \dotsc, n}$ such that $b_{i_0} \neq 0$. Now a function $\Psi$ with the desired property is for example given by
\bgloz
  \Psi(\sum_{i=i_0}^n b_i \omega_i + T R\sig) \defeq \psi(b_{i_0}).
\egloz
\eproof
This completes the proof of Lemma~\ref{i_*}.
\eproof
Using
\bgloz
  C(\overline{R\sig}) \rtimes R\sig \rtimes \Fz_q\reg \cong \ilim_m \gekl{C(R\sig / T^m R\sig) \rtimes R_m\sig \rtimes \Fz_q\reg; i_{m,m+1}}
\egloz
and continuity of $K_0$, the previous lemma yields
\blemma
\label{K(R&F)}
We can identify $K_0(C(\overline{R\sig}) \rtimes R\sig \rtimes \Fz_q\reg)$ with $\Zz[\tfrac{1}{q}] \oplus \bigoplus_{\dualF \setminus \gekl{1}} \Zz$ such that for every $m \in \Zz \nneg$ and $\chi \in \dualF \setminus \gekl{1}$, $[\1z_m]$ corresponds to 
$
  \rukl{
  \begin{smallmatrix}
  q^{-mn} \\
  0 \\
  \vdots \\
  0
  \end{smallmatrix}
  } 
$
and $[\1z_m p_\chi]$ corresponds to
$
  \rukl{
  \begin{smallmatrix}
  q^{-mn} Q^{mn}\\
  0 \\
  \vdots \\
  0 \\
  1 \\
  0 \\
  \vdots \\
  0
  \end{smallmatrix}
  }
$.
Here, as in Proposition~\ref{K(add&rou)}, the entry \an{1} in the vector corresponding to $[\1z_m p_\chi]$ is the $\chi$-th entry in $\bigoplus_{\dualF \setminus \gekl{1}} \Zz$, and $Q^{mn} = \tfrac{q^{mn}-1}{q-1}$.
\elemma
Now, it is only a small step from this lemma to our proposition. We just have to invert $\Fz_q[T]\reg = \Fz_q[T] \setminus \gekl{0}$ in $R\sig$. Namely, we know that
\bgloz
  \Az_T\sig = (\Fz_q[T]\reg)^{-1} \overline{R\sig} \text{ and } K\sig = (\Fz_q[T]\reg)^{-1} R\sig.
\egloz
This implies
\bgloz
  C_0(\Az_T\sig) \rtimes K\sig \rtimes \Fz_q\reg \cong \ilim_{\Fz_q[T]\reg} \gekl{C(\overline{R\sig}) \rtimes R\sig \rtimes \Fz_q\reg; \mu}
\egloz
where $\mu$ is the multiplicative action of the semigroup $\Fz_q[T]\reg$ on $C(\overline{R\sig}) \rtimes R\sig \rtimes \Fz_q\reg$ via endomorphisms. This observation  already implies, using continuity of $K_1$, that $K_1(C_0(\Az_T\sig) \rtimes K\sig \rtimes \Fz_q\reg)$ vanishes since we know that $K_1(C(\overline{R\sig}) \rtimes R\sig \rtimes \Fz_q\reg)$ vanishes.

Moreover, it turns out that for every $f$ in $\Fz_q[T]\reg$, $\mu_f$ is an isomorphism on $K_0(C(\overline{R\sig}) \rtimes R\sig \rtimes \Fz_q\reg)$. Namely, with the help of Lemma~\ref{K(R&F)}, it is easy to compute that
\bgloz
  (\mu_f)_* = \id \text{ on } K_0(C(\overline{R\sig}) \rtimes R\sig \rtimes \Fz_q\reg) \text{ for every } f \in \Fz_q[T]\reg \text{ with } v_T(f) = 0
\egloz
and
\bgloz
  (\mu_T)_* = 
  \rukl{
  \begin{smallmatrix}
  q^{-n} & -q^{-n}Q^n & \dotso & -q^{-n}Q^n \\
  0 & 1 & & 0 \\
  \vdots & & \ddots & \\
  0 & 0 & & 1
  \end{smallmatrix}
  }  
\egloz
with respect to the identification $K_0(C(\overline{R\sig}) \rtimes R\sig \rtimes \Fz_q\reg) \cong \Zz[\tfrac{1}{q}] \oplus \bigoplus_{\dualF \setminus \gekl{1}} \Zz$ from Lemma~\ref{K(R&F)}. Again, $Q^n = \tfrac{q^n-1}{q-1}$.

But this implies, again by continuity of $K_0$, that the canonical homomorphism
\bgloz
  C(\overline{R\sig}) \rtimes R\sig \rtimes \Fz_q\reg \to C_0(\Az_T\sig) \rtimes K\sig \rtimes \Fz_q\reg
\egloz
induces an isomorphism on $K_0$.

Proposition~\ref{K(add&rou)} follows immediately from this observation and Lemma~\ref{K(R&F)}.

\section{Adjoining a prime element}
\label{prime}
To pass from $C_0(\Az_T\sig) \rtimes K\sig \rtimes \Fz_q\reg$ to $C_0(\Az_T\sig) \rtimes K\sig \rtimes (K\sig)\reg$, we have to treat the multiplicative action, say $\mu$. We proceed in two steps which occupy the next two sections.

From now on, we assume that $K$ has only one single infinite place. This means that $K\sig$ has only one single place $w$ sitting above $v_T$.

Our first goal is to treat the multiplicative action of a prime element for $w$. Take such a prime element in $R\sig$, say $\tau$. We observe that we can always find a free abelian subgroup $\Gamma'$ of $(K\sig)\reg$ such that $(K\sig)\reg = \Fz_q\reg \times \spkl{\tau} \times \Gamma'$. Later on, we will choose a particular $\Gamma'$ with special properties.

Let us now determine $(\mu_\tau)_*$ on $K_0(C_0(\Az_T\sig) \rtimes K\sig \rtimes \Fz_q\reg)$. Here is the result:

Let $f$ be the inertia degree at $w$, i.e. $(R\sig : \tau R\sig) = q^f$.
\bprop
\label{mu_psi}
With respect to the identification 
\bgloz
  K_0(C_0(\Az_T\sig) \rtimes K\sig \rtimes \Fz_q\reg) \cong \Zz[\tfrac{1}{q}] \oplus \bigoplus_{\dualF \setminus \gekl{1}} \Zz
\egloz
from Proposition~\ref{K(add&rou)}, $(\mu_\tau)_*$ is of the form
\bgloz
  (\mu_\tau)_* = 
  \rukl{
  \begin{smallmatrix}
  q^{-f} & -q^{-f}Q^f & \dotso & -q^{-f}Q^f \\
  0 & 1 & & 0 \\
  \vdots & & \ddots & \\
  0 & 0 & & 1
  \end{smallmatrix}
  }.
\egloz
\eprop
\bproof
Proposition~\ref{K(add&rou)} tells us that the copy of $\Zz[\tfrac{1}{q}]$ is generated by the $K_0$-classes $\menge{\eckl{\1z_m}}{m \in \Zz \nneg}$ and that the $q-2$ copies of $\Zz$ (i.e. $\bigoplus_{\dualF \setminus \gekl{1}} \Zz$) are generated by the $K_0$-classes $\menge{\eckl{\1z p_\chi}}{1 \neq \chi \in \dualF}$ with respect to the identification
\bgloz
  K_0(C_0(\Az_T\sig) \rtimes K\sig \rtimes \Fz_q\reg) \cong \Zz[\tfrac{1}{q}] \oplus \bigoplus_{\dualF \setminus \gekl{1}} \Zz
\egloz
from Proposition~\ref{K(add&rou)}.

We have
\bgloz
  \mu_\tau(\1z_m) = \1z_{\tau T^m \overline{R\sig}} = \sum_{b \in \tau R\sig / T R\sig} v^{T^m b} \1z_{m+1} v^{-T^m b}.
\egloz
As each of the pairwise orthogonal projections $v^{T^m b} \1z_{m+1} v^{-T^m b}$ gives the same class as $\1z_{m+1}$ in $K_0$, we obtain
\bglnoz
  (\mu_\tau)_*(\eckl{\1z_m}) &=& (\tau R\sig : T R\sig) \eckl{\1z_{m+1}} \\
  &=& (R\sig : T R\sig) (R\sig : \tau R\sig)^{-1} \eckl{\1z_{m+1}} \\
  &=& q^n q^{-f} \eckl{\1z_{m+1}} = q^{n-f} \eckl{\1z_{m+1}}.
\eglnoz
As $\eckl{\1z_{m}} \hat{=} q^{-mn}$ and $\eckl{\1z_{m+1}} \hat{=} q^{-(m+1)n}$, we obtain that $(\mu_\tau)_*$ acts by multiplication with $q^{-f}$ on the copy of $\Zz[\tfrac{1}{q}]$.

It remains to determine $(\mu_\tau)_*(\eckl{\1z p_\chi})$ for $1 \neq \chi \in \dualF$. To do so, we construct for every $\psi$, $\chi$ in $\dualF$ a partial isometry with support projection $(\1z_{\tau \overline{R\sig}} - \1z_{T \overline{R\sig}})p_{\overline{\psi}\chi}$ and range projection $(\1z_{\tau \overline{R\sig}} - \1z_{T \overline{R\sig}})p_\chi$. The procedure is very similar to the one in the proof of Lemma~\ref{i_*}. What we need is an element $x_\Psi \in C_0(\Az_T\sig)$ corresponding to $\psi \in \dualF$ with the properties
\bgloz
  x_\Psi x_\Psi^* = x_\Psi^* x_\Psi = \1z_{\tau \overline{R\sig}} - \1z_{T \overline{R\sig}}
\egloz
and
\bgloz
  p_\chi x_\Psi = x_\Psi p_{\overline{\psi}\chi} \fa \chi \in \dualF.
\egloz
Let $(\tau R\sig / T R\sig)\reg$ be $(\tau R\sig / T R\sig) \setminus \gekl{0 + T R\sig}$. A similar computation as in Lemma~\ref{i_*} shows that an element $x_\Psi$ with the desired properties is given by
\bgloz
  x_\Psi = \sum_{b + T R\sig \in (\tau R\sig / T R\sig)\reg} \Psi(b + T R\sig) v^b \1z_{T \overline{R\sig}} v^{-b}
\egloz
provided that $\Psi$ is a function $(\tau R\sig / T R\sig)\reg \to \Tz$ satisfying
\bgloz
  \Psi(a \cdot b + T R\sig) = \psi(a) \Psi(b + T R\sig) \fa a \in \Fz_q\reg.
\egloz
Such a function $\Psi$ exists by Lemma~\ref{Psi-exists}: We just have to restrict the function constructed in that lemma to
\bgloz
  (\tau R\sig / T R\sig)\reg \subseteq (R\sig / T R\sig)\reg.
\egloz
Thus, there is a partial isometry with the desired properties. This implies
\bgl
\label{pp'}
  \eckl{(\1z_{\tau \overline{R\sig}} - \1z_{T \overline{R\sig}})p_{\overline{\psi}\chi}} 
  = \eckl{(\1z_{\tau \overline{R\sig}} - \1z_{T \overline{R\sig}})p_\chi}
  \fa \psi, \chi \in \dualF.
\egl
We can then compute as in the proof of Lemma~\ref{i_*}:
\bgl
\label{mu(p)_1}
  \eckl{\mu_\tau(\1z p_\chi)} = \eckl{\1z_{\tau \overline{R\sig}} p_\chi} 
  = \eckl{(\1z_{\tau \overline{R\sig}} - \1z_{T \overline{R\sig}})p_\chi} + \eckl{\1z_{T \overline{R\sig} p_\chi}}
\egl
and
\bgln
\label{mu(p)_2}
  && (q-1)\eckl{(\1z_{\tau \overline{R\sig}}-\1z_{T \overline{R\sig}})p_\chi} \\
  &\overset{\eqref{pp'}}{=}& \sum_{\psi \in \dualF} \eckl{(\1z_{\tau \overline{R\sig}}-\1z_{T \overline{R\sig}})p_{\overline{\psi}\chi}} \nonumber \\
  &=& \eckl{\1z_{\tau \overline{R\sig}}-\1z_{T \overline{R\sig}}} = (q^{n-f}-1)\eckl{\1z_{T \overline{R\sig}}} \nonumber.
\egln
\eqref{mu(p)_2} implies that
\bgl
\label{mu(p)_3}
  \eckl{(\1z_{\tau \overline{R\sig}}-\1z_{T \overline{R\sig}})p_\chi} = Q^{n-f} \eckl{\1z_{T \overline{R\sig}}} \fa \chi \in \dualF.
\egl
Inserting \eqref{mu(p)_3} in \eqref{mu(p)_1} yields
\bgloz
  \eckl{\mu_\tau(\1z p_\chi)} = Q^{n-f} \eckl{\1z_{T \overline{R\sig}}} + \eckl{\1z_{T \overline{R\sig}}p_\chi}.
\egloz
In the identification
\bgloz
  K_0(C_0(\Az_T\sig) \rtimes K\sig \rtimes \Fz_q\reg) \cong \Zz[\tfrac{1}{q}] \oplus \bigoplus_{\dualF \setminus \gekl{1}} \Zz
\egloz
from Proposition~\ref{K(add&rou)}, $\eckl{\1z_{T \overline{R\sig}}} = \eckl{\1z_1}$ corresponds to
$
  \rukl{
  \begin{smallmatrix}
  q^{-n} \\
  0 \\
  \vdots \\
  0
  \end{smallmatrix}
  } 
$
and $\eckl{\1z_{T \overline{R\sig}}p_\chi} = \eckl{\1z_1 p_\chi}$ corresponds to
$
  \rukl{
  \begin{smallmatrix}
  q^{-n} Q^n \\
  0 \\
  \vdots \\
  0 \\
  1 \\
  0 \\
  \vdots \\
  0
  \end{smallmatrix}
  }
$.
We therefore obtain
\bgloz
  \eckl{\mu_\tau(\1z p_\chi)} \hat{=}
  \rukl{
  \begin{smallmatrix}
  q^{-n} Q^{n-f} - q^{-n} Q^n \\
  0 \\
  \vdots \\
  0 \\
  1 \\
  0 \\
  \vdots \\
  0
  \end{smallmatrix}
  }
  =
  \rukl{
  \begin{smallmatrix}
  \tfrac{1}{q^n} \tfrac{q^{n-f}-1-(q^n-1)}{q-1} \\
  0 \\
  \vdots \\
  0 \\
  1 \\
  0 \\
  \vdots \\
  0
  \end{smallmatrix}
  }
  =
  \rukl{
  \begin{smallmatrix}
  -q^{-f} Q^f \\
  0 \\
  \vdots \\
  0 \\
  1 \\
  0 \\
  \vdots \\
  0
  \end{smallmatrix}
  }.
\egloz
\eproof

Going through the Pimsner-Voiculescu sequence, we can now derive the following
\bprop
\label{K(AKFtau)}
We have
\bgloz
  K_0(C_0(\Az_T\sig) \rtimes K\sig \rtimes (\Fz_q\reg \times \spkl{\tau})) \cong (\Zz / Q^f \Zz) \oplus \bigoplus_{\chi \in \dualF \setminus \gekl{1}} \Zz
\egloz
and
\bgloz
  K_1(C_0(\Az_T\sig) \rtimes K\sig \rtimes (\Fz_q\reg \times \spkl{\tau})) \cong \bigoplus_{\chi \in \dualF \setminus \gekl{1}} \Zz.
\egloz
Moreover, the $K_0$-class $\eckl{\1z}$ is a generator for the copy of $\Zz / Q^f \Zz$, and generators in $K_0$ of $\bigoplus_{\chi \in \dualF \setminus \gekl{1}} \Zz$ are given by the classes $\menge{\eckl{\1z p_\chi}}{\chi \in \dualF \setminus \gekl{1}}$. Generators for $K_1$ are given by the classes $\menge{\eckl{w_\chi}}{\chi \in \dualF \setminus \gekl{1}}$ where $w_\chi$ are unitaries in the unitalization of $C_0(\Az_T\sig) \rtimes K\sig \rtimes (\Fz_q\reg \times \spkl{\tau})$ of the form
\bgl
\label{w}
  w_\chi = t_\tau (\1z p_1) + x_X p_\chi + (\1z p_\chi)t_\tau^* + (1-\1z p_1-\1z p_\chi).
\egl
Here, we can take for $x_X$ any element in $C_0(\Az_T\sig)$ satisfying
\bgloz
  x_X x_X^* = x_X^* x_X = \1z - \1z_{\tau \overline{R\sig}}
\egloz
and
\bgloz
  x_X p_\chi = p_1 x_X.
\egloz
\eprop
The proof of this proposition is analogous to the proof of Proposition~5.10 in \cite{Cu-Li3}.

\bremark
\label{x_X}
The element $x_X$ can be chosen as
\bgloz
  x_X = \sum_{b + \tau R\sig \in (R\sig / \tau R\sig)\reg} X(b + \tau R\sig) v^b \1z_{\tau \overline{R\sig}} v^{-b}
\egloz
where $X$ is a function $X: (R\sig / \tau R\sig)\reg \to \Tz$ with the property that
\bgloz
  X(a \cdot b + \tau R\sig) = \overline{\chi}(a) X(b + \tau R) \fa a \in \Fz_q\reg.
\egloz
Such a function exists as $R\sig / \tau R\sig \cong \Fz_{q^f}$ and we can identify $\Fz_{q^f}$ with $(\Fz_q)^f$ as $\Fz_q$-vector spaces. Then we just have to apply the analogous construction as in the proof of Lemma~\ref{Psi-exists}.

Later on, we will choose particular elements $x_X$ satisfying an invariance property.
\eremark

In the case of polynomial rings and their quotient fields, it was possible to go on from this point and to compute the K-theory of the associated ring C*-algebras. The idea was to compare the situation with higher dimensional tori and to use the Pimsner-Voiculescu sequence. Since there was no torsion, everything could be computed.

Now, as we have seen, torsion can appear. This leads do extension problems which we cannot solve up to now. Another complication that arises in our more general situation is that the field of constants $\Fz_q$ in $K$ (or $K\sig$) does not need to coincide with the field of constants of the local field at infinity (or at the place of $K\sig$ over $v_T$). Therefore, it is not clear any more whether we can choose generators for the remaining part of the multiplicative group which fix the unitaries $w_\chi$ (at the level of C*-algebras, not just in K-theory). This step was however crucial in our computations for polynomial rings. Note that both complications arise if and only if the inertia degree at infinity of our function field $K / \Fz_q(T)$ (or of $K\sig / \Fz_q(T)$ over $v_T$) is strictly bigger than $1$.

In the sequel, we proceed in two directions: First, it turns out that once we rationalize, these two problems disappear. This is of course not surprising for the extension problems caused by torsion, but the second problem can be solved as well. And secondly, we develop a solution of the second problem under a certain primeness condition, namely that the inertia degree $f$ at infinity of $K / \Fz_q(T)$ (or of $K\sig / \Fz_q(T)$ over $v_T$) satisfies $(f,q-1)=1$. In this case, we can say more about torsion, and this leads to our observation that the K-theory of ring C*-algebras detects the inertia degrees at infinity.

\section{Rational K-theory computations}
\label{rational}
Our goal is to prove Theorem~\ref{mainthm1}. Here is the statement:

Take a function field $K / \Fz_p(T)$ with fields of constants $\Fz_q$. Assume that the field extension $K / \Fz_q(T)$ has only one single infinite place with inertia degree $f$. Define $Q^f = \tfrac{q^f-1}{q-1}$. Let $\Gamma \subseteq K\reg$ be a free abelian subgroup such that $K\reg = \Fz_q\reg \times \Gamma$.
\btheo
After inverting $Q^f$, the K-theory of $\fA[R]$, the ring C*-algebra associated to the ring of integers $R$ in $K / \Fz_q(T)$, is given by
\bgloz
  K_*(\fA[R]) \otimes_{\Zz} \Zz[\tfrac{1}{Q^f}] \cong \rukl{\ti{K}_0(C^*(\Fz_q\reg)) \otimes_{\Zz} \extalg(\Gamma)} \otimes_{\Zz} \Zz[\tfrac{1}{Q^f}].
\egloz
\etheo

Let us come to the proof of this theorem. Since we have already computed K-theory for $C_0(\Az_T\sig) \rtimes K\sig \rtimes (\Fz_q\reg \times \spkl{\tau})$, it remains to treat the multiplicative action of a subgroup $\Gamma'$ of $(K\sig)\reg$ with the property that $(K\sig)\reg = \Fz_q\reg \times \spkl{\tau} \times \Gamma'$. Now, let us fix such a subgroup $\Gamma'$. $\mu$ denotes the multiplicative action of $\Gamma'$ on $C_0(\Az_T\sig) \rtimes K\sig \rtimes (\Fz_q\reg \times \spkl{\tau})$. We start with the following
\blemma
\label{mu_g=id}
Let $\gamma_1, \dotsc, \gamma_i$ and $\gamma$ be linearly independent elements in $\Gamma'$. Assume that there are unitaries $\menge{w_\chi}{\chi \in \dualF \setminus \gekl{1}}$ as in Proposition~\ref{K(AKFtau)} which are fixed by $\mu_{\gamma_1}, \dotsc, \mu_{\gamma_i}$ and $\mu_\gamma$. Then we have $(\mu_\gamma)_* = \id$ on $K_*(C_0(\Az_T\sig) \rtimes K\sig \rtimes (\Fz_q\reg \times (\spkl{\tau} \times \spkl{\gamma_1, \dotsc, \gamma_i}))$.
\elemma
\bproof
Consider the Pimsner-Voiculescu sequence for
\bglnoz
  && C_0(\Az_T\sig) \rtimes K\sig \rtimes (\Fz_q\reg \times (\spkl{\tau} \times \spkl{\gamma_1, \dotsc, \gamma_i} \times \spkl{\gamma})) \\
  &\cong& C_0(\Az_T\sig) \rtimes K\sig \rtimes (\Fz_q\reg \times (\spkl{\tau} \times \spkl{\gamma_1, \dotsc, \gamma_i})) \rtimes_{\mu_\gamma} \Zz.
\eglnoz
Exactness of the Pimsner-Voiculescu sequence implies that our assertion is equivalent to saying that the boundary maps are surjective. To show this, we construct K-theoretic elements
\bgloz
  \eckl{\1z, t_{\gamma_{j_1}}, t_{\gamma_{j_2}}, \dotsc, t_{\gamma_{j_l}}}, \eckl{\1z \cdot p_\chi, t_{\gamma_{j_1}}, t_{\gamma_{j_2}}, \dotsc, t_{\gamma_{j_l}}}, 
  \eckl{w_\chi, t_{\gamma_{j_1}}, t_{\gamma_{j_2}}, \dotsc, t_{\gamma_{j_l}}}
\egloz
and
\bgloz
  \eckl{\1z, t_{\gamma_{j_1}}, \dotsc, t_{\gamma_{j_l}}, t_{\gamma}}, \eckl{\1z \cdot p_\chi, t_{\gamma_{j_1}}, \dotsc, t_{\gamma_{j_l}}, t_{\gamma}}, 
  \eckl{w_\chi, t_{\gamma_{j_1}}, \dotsc, t_{\gamma_{j_l}}, t_{\gamma}}
\egloz
for every $1 \leq j_i < \dotsb < j_l \leq i$, $0 \leq j$ and $\chi \in \dualF \setminus \gekl{1}$ by comparing our situation with higher dimensional tori as in \cite{Cu-Li3}, Section~6. This is possible because of our assumption that the elements $w_\chi$ are invariant under $\mu_{\gamma_1}, \dotsc, \mu_{\gamma_i}$ and $\mu_{\gamma}$. Now the claim is that
\bglnoz
  && \menge{\eckl{\1z, t_{\gamma_{j_1}}, t_{\gamma_{j_2}}, \dotsc, t_{\gamma_{j_l}}}}{1 \leq j_1 < \dotsb < j_l \leq i, 0 \leq l} \\
  &\cup& \menge{\eckl{\1z \cdot p_\chi, t_{\gamma_{j_1}}, t_{\gamma_{j_2}}, \dotsc, t_{\gamma_{j_l}}}}
  {1 \leq j_1 < \dotsb < j_l \leq i, 0 \leq l; \chi \in \dualF \setminus \gekl{1}}
  \\
  &\cup& \menge{\eckl{w_\chi, t_{\gamma_{j_1}}, t_{\gamma_{j_2}}, \dotsc, t_{\gamma_{j_l}}}}
  {1 \leq j_1 < \dotsb < j_l \leq i, 0 \leq l; \chi \in \dualF \setminus \gekl{1}}
\eglnoz
are generators for $K_*(C_0(\Az_T\sig) \rtimes K\sig \rtimes (\Fz_q\reg \times (\spkl{\tau} \times \spkl{\gamma_1, \dotsc, \gamma_i}))$ and that the boundary maps in the Pimsner-Voiculescu sequence for
\bglnoz
  && C_0(\Az_T\sig) \rtimes K\sig \rtimes (\Fz_q\reg \times (\spkl{\tau} \times \spkl{\gamma_1, \dotsc, \gamma_i} \times \spkl{\gamma})) \\
  &\cong& C_0(\Az_T\sig) \rtimes K\sig \rtimes (\Fz_q\reg \times (\spkl{\tau} \times \spkl{\gamma_1, \dotsc, \gamma_i})) \rtimes_{\mu_\gamma} \Zz
\eglnoz
send
\bglnoz
  && \eckl{\1z, t_{\gamma_{j_1}}, t_{\gamma_{j_2}}, \dotsc, t_{\gamma_{j_l}}, t_\gamma} 
  \text{ to } \eckl{\1z, t_{\gamma_{j_1}}, t_{\gamma_{j_2}}, \dotsc, t_{\gamma_{j_l}}}, \\
  && \eckl{\1z \cdot p_\chi, t_{\gamma_{j_1}}, t_{\gamma_{j_2}}, \dotsc, t_{\gamma_{j_l}}, t_\gamma} 
  \text{ to } \eckl{\1z \cdot p_\chi, t_{\gamma_{j_1}}, t_{\gamma_{j_2}}, \dotsc, t_{\gamma_{j_l}}} \\
  &\text{and}& \eckl{w_\chi, t_{\gamma_{j_1}}, t_{\gamma_{j_2}}, \dotsc, t_{\gamma_{j_l}}, t_\gamma} 
  \text{ to } \eckl{w_\chi, t_{\gamma_{j_1}}, t_{\gamma_{j_2}}, \dotsc, t_{\gamma_{j_l}}}
\eglnoz
for every $1 \leq j_1 < \dotsb < j_l \leq i, 0 \leq l; \chi \in \dualF \setminus \gekl{1}$.

It actually suffices to prove the second assertion because the first one can be deduced inductively as follows: The starting point is Proposition~\ref{K(AKFtau)}, and the induction step from $j$ to $j+1$ is given by the second observation applied to $\gamma_1, \dotsc, \gamma_j$ and $\gamma = \gamma_{j+1}$. 

The proof of the second assertion is the same as the one for Proposition~7.1 in \cite{Cu-Li3}.
\eproof
Now, we can choose $\Gamma' \subseteq (K\sig)\reg$ so that ${(K\sig)}^{\gekl{w}} = \Fz_q \times \Gamma'$ where ${K\sig}^{\gekl{w}} = \menge{a \in K\sig}{w(a)=0}$. By construction, the inclusion $K\sig \into K_w\sig = \Az_T\sig$ induces a map $\Gamma' \to R_w\sig = \overline{R\sig} \to \overline{R\sig} / \tau \overline{R\sig}$ with image in $(\overline{R\sig} / \tau \overline{R\sig})\reg$, i.e. a group homomorphism $\Gamma' \to (\overline{R\sig} / \tau \overline{R\sig})\reg; \gamma \ma \overline{\gamma}$. It is clear that if we choose $x_X$ in the construction of $w_\chi$ to be of the form as in Remark~\ref{x_X}, the unitaries $w_\chi$ for all $\chi \in \dualF \setminus \gekl{1}$ will be fixed by $\mu_\gamma$ if $\overline{\gamma} = 1 + \tau \overline{R\sig}$. So, let us choose particular generators $\gamma_1, \gamma_2, \dotsc$ for a subgroup $\Gamma' \subseteq (K\sig)\reg$ with ${(K\sig)}^{\gekl{w}} = \Fz_q \times \Gamma'$. We can arrange that
\bgloz
  \overline{\gamma_1} \text{ is a multiplicative generator of } (\overline{R\sig} / \tau \overline{R\sig})\reg
\egloz
and that
\bgloz
  \overline{\gamma_i} = 1 + \tau \overline{R\sig} \fa i \geq 2.
\egloz
Set $\Gamma_i \defeq \spkl{\tau, \gamma_1, \dotsc, \gamma_{i-1}}$.
\bprop
We have
\bgloz
  K_*(C_0(\Az_T\sig) \rtimes K\sig \rtimes (\Fz_q\reg \times \Gamma_i)) \otimes_{\Zz} \Zz[\tfrac{1}{Q^f}] \cong 
  \rukl{\ti{K}_0(C^*(\Fz_q\reg)) \otimes_{\Zz} \extalg(\Gamma_i)} \otimes_{\Zz} \Zz[\tfrac{1}{Q^f}]
\egloz
for all $i \geq 2$ and under these identifications, the canonical map 
\bgloz
  C_0(\Az_T\sig) \rtimes K\sig \rtimes (\Fz_q\reg \times \Gamma_i) \to C_0(\Az_T\sig) \rtimes K\sig \rtimes (\Fz_q\reg \times \Gamma_{i+1})
\egloz
induces the canonical inclusion
\bgloz
  \rukl{\ti{K}_0(C^*(\Fz_q\reg)) \otimes_{\Zz} \extalg(\Gamma_i)} \otimes_{\Zz} \Zz[\tfrac{1}{Q^f}] 
  \to \rukl{\ti{K}_0(C^*(\Fz_q\reg)) \otimes_{\Zz} \extalg(\Gamma_{i+1})} \otimes_{\Zz} \Zz[\tfrac{1}{Q^f}]
\egloz
for every $i \geq 2$.
\eprop
\bproof
By the previous lemma and the Pimsner-Voiculescu sequence, we obtain
\bglnoz
  && K_*(C_0(\Az_T\sig) \rtimes K\sig \rtimes (\Fz_q\reg \times \spkl{\tau, \gamma_2, \dotsc, \gamma_{i-1}})) \otimes_{\Zz} \Zz[\tfrac{1}{Q^f}] \\
  &\cong& 
  \rukl{\ti{K}_0(C^*(\Fz_q\reg)) \otimes_{\Zz} \extalg(\spkl{\tau, \gamma_2, \dotsc, \gamma_{i-1}})} \otimes_{\Zz} \Zz[\tfrac{1}{Q^f}] \fa i \geq 2
\eglnoz
because $\overline{\gamma_i} = 1 + \tau \overline{R\sig}$ for all $i \geq 2$ implies $\mu_{\gamma_i}(w_\chi)=w_\chi$ for all $\chi \in \dualF \setminus \gekl{1}$ and $i \geq 2$ as observed above.

It remains to prove that $(\mu_{\gamma_1})_* \otimes \id_{\Zz[\tfrac{1}{Q^f}]} = \id \otimes \id_{\Zz[\tfrac{1}{Q^f}]}$ on $K_*(C_0(\Az_T\sig) \rtimes K\sig \rtimes (\Fz_q\reg \times \spkl{\tau, \gamma_2, \dotsc, \gamma_{i-1}})) \otimes_{\Zz} \Zz[\tfrac{1}{Q^f}]$. We know that $\overline{\gamma_1} \in (\overline{R\sig} / \tau \overline{R\sig})\reg$ has order $q^f-1$ because $(\overline{R\sig} / \tau \overline{R\sig})\reg \cong \Fz_{q^f} \reg$ is a cyclic group of order $q^f-1$. Thus $\overline{(\gamma_1)^{q^f-1}} = 1 + \tau \overline{R\sig}$, and this implies $\mu_{(\gamma_1)^{q^f-1}}(w_\chi) = w_\chi$ for all $\chi \in \dualF \setminus \gekl{1}$ as explained above. Hence, by our previous lemma, we have
\bgloz
  (\mu_{(\gamma_1)^{q^f-1}})_* \otimes \id_{\Zz[\tfrac{1}{Q^f}]} = \id \otimes \id_{\Zz[\tfrac{1}{Q^f}]}
\egloz
on $K_*(C_0(\Az_T\sig) \rtimes K\sig \rtimes (\Fz_q\reg \times \spkl{\tau, \gamma_2, \dotsc, \gamma_{i-1}})) \otimes_{\Zz} \Zz[\tfrac{1}{Q^f}]$. Thus
\bgloz
  \rukl{(\mu_{\gamma_1})_* \otimes \id_{\Zz[\tfrac{1}{Q^f}]}}^{q^f-1} 
  = \rukl{(\mu_{\gamma_1})_*}^{q^f-1} \otimes \id_{\Zz[\tfrac{1}{Q^f}]} 
  = \id \otimes \id_{\Zz[\tfrac{1}{Q^f}]}.
\egloz
This tells us that $(\mu_{\gamma_1})_* \otimes \id_{\Zz[\tfrac{1}{Q^f}]}$ is of order $q^f-1$, and with the help of the Pimsner-Voiculescu sequence (proceeding inductively on $i$), this implies $(\mu_{\gamma_1})_* \otimes \id_{\Zz[\tfrac{1}{Q^f}]} = \id \otimes \id_{\Zz[\tfrac{1}{Q^f}]}$ on $K_*(C_0(\Az_T\sig) \rtimes K\sig \rtimes (\Fz_q\reg \times \spkl{\tau, \gamma_2, \dotsc, \gamma_{i-1}})) \otimes_{\Zz} \Zz[\tfrac{1}{Q^f}]$ for all $i \geq 2$.
\eproof
Now, let us set $\Gamma \defeq \spkl{\tau} \times \Gamma'$ so that $K\reg = (K\sig)\reg = \Fz_q\reg \times \Gamma$. By continuity of K-theory, we deduce from the previous proposition
\bgloz
  K_*(C_0(\Az_T\sig) \rtimes K\sig \rtimes (K\sig)\reg) \otimes_{\Zz} \Zz[\tfrac{1}{Q^f}] \cong 
  \rukl{\ti{K}_0(C^*(\Fz_q\reg)) \otimes_{\Zz} \extalg(\Gamma)} \otimes_{\Zz} \Zz[\tfrac{1}{Q^f}].
\egloz
This result, together with \eqref{A-Mor-infadele} and \eqref{AKK-Mor-AKsigKsig}, completes the proof of Theorem~\ref{mainthm1}.

Note that Theorem~\ref{mainthm1} covers the results in \cite{Cu-Li3} since for the function fields $\Fz_q(T) / \Fz_p(T)$, the inertia degree at infinity $f$ is $1$, i.e. $Q^f=1$.

An immediate consequence of Theorem~\ref{mainthm1} is
\bcor
\label{torsion-order}
In the same situation as in Theorem~\ref{mainthm1}, let $x \in K_*(\fA[R])$ be a torsion element of order $N$. Then $N$ divides a suitable power of $Q^f$, i.e. $N \tei (Q^f)^j$ for $j$ big enough. In particular, every prime divisor of $N$ is a prime divisor of $Q^f$.
\ecor
\bremark
It becomes clear why we restrict our investigations to the case of one single infinite place. The reason is that if we allow more infinite places, we will have to deal with many prime elements, one for each infinite place. But then, it is not possible to carry out the explicit computation in Section~\ref{prime} for the crossed product involving the action of all these prime elements. However, the explicit description of K-theory at this stage is the basis for what we do in this section, so it is not clear how to proceed in the case of more than one infinite place.
\eremark
Of course, it would be desirable to say more about the torsion part of $K_*(\fA[R])$. We do so in the next section, but even there we are far away from completely determining the torsion part. So this remains an open problem.

\section{Torsion in K-theory determines inertia degree}
\label{torsion}
In this section, we consider the same situation as in the previous section. In addition, we assume that the inertia degree $f$ at infinity of our function field $K / \Fz_q(T)$ is prime to the cardinality of the multiplicative group of the field of constants $\Fz_q$, i.e. $(f,q-1)=1$. In this case, we will be able to choose a particular element $x_X$ or rather a particular function $X$ such that the corresponding unitaries $w_\chi$ are fixed by every element in $\Gamma'$ for a suitable choice of $\Gamma'$. The unitaries $w_\chi$ were defined in Proposition~\ref{K(AKFtau)}, see \eqref{w}, and the elements $x_X$ were constructed in Remark~\ref{x_X}.

First of all, note that 
\bgloz
  Q^f = \tfrac{q^f-1}{q-1} = 1+q+\dotsb+q^{f-1} \equiv f \mod q-1.
\egloz
Hence it follows from $(f,q-1)=1$ that we also have $(Q^f,q-1)=1$. Thus the homomorphism $\Fz_q\reg \to \Fz_q\reg; a \ma a^{Q^f}$ is an isomorphism.

Now let $\chi$ be a character in $\dualF$. Composing $\chi$ by the inverse of the map $\Fz_q\reg \to \Fz_q\reg; a \ma a^{Q^f}$, we obtain a character $\ti{\chi}$ with the property $\ti{\chi}(a^{Q^f}) = \chi(a)$ for every $a$ in $\Fz_q\reg$. Since $Q^f \cdot (q-1) = q^f-1$ is the order of the cyclic group $(R\sig/ \tau R\sig)\reg$, we know that $(b + \tau R\sig)^{Q^f}$ lies in $\Fz_q\reg \subseteq (R\sig/ \tau R\sig)\reg$ for every element $b + \tau R\sig \in (R\sig/ \tau R\sig)\reg$. Thus we can define $X: (R\sig/ \tau R\sig)\reg \to \Tz$ as the composition $(R\sig/ \tau R\sig)\reg \overset{(\cdot)^{Q^f}}{\to} \Fz_q\reg \overset{\ti{\chi}}{\to} \Tz$, i.e.
\bgloz
  X(b + \tau R\sig) = \ti{\chi}((b + \tau R\sig)^{Q^f}).
\egloz
By construction, we have
\bgloz
  X(a \cdot b + \tau R\sig) = \ti{\chi}(a^{Q^f}) \ti{\chi}((b + \tau R\sig)^{Q^f}) = \chi(a) X(b + \tau R\sig).
\egloz
The point is that we can now choose $\Gamma'$ with particular free generators so that the unitaries $w_\chi$ are fixed by the multiplicative action of $\Gamma'$ for every $\chi \in \dualF \setminus \gekl{1}$. Recall that $\Gamma'$ is defined as a subgroup of $(K\sig)\reg$ with $(K\sig)\reg = \Fz_q\reg \times \spkl{\tau} \times \Gamma'$ or, equivalently, ${(K\sig)}^{\gekl{w}} = \Fz_q \times \Gamma'$. Let us start with an arbitrary choice of $\Gamma'$ together with free generators $\gekl{\gamma_i}$ for $\Gamma'$. We know that $(\gamma_i + \tau R\sig)^{Q^f}$ lies in $\Fz_q\reg \subseteq (R\sig/ \tau R\sig)\reg$, so there exists an element $a_i$ in $\Fz_q\reg$ such that
\bgl
\label{ag^Q}
  (a_i \gamma_i + \tau R\sig)^{Q^f} = 1 + \tau R\sig \text{ in } (R\sig/ \tau R\sig)\reg.
\egl
(Recall that $a \ma a^{Q^f}$ is an isomorphism of $\Fz_q\reg$.) Thus we obtain
\bglnoz
  \mu_{a_i \gamma_i}(x_X) &=& \sum_{b + \tau R\sig \in (R\sig / \tau R\sig)\reg} X(b + \tau R\sig) v^{a_i \gamma_i b} \1z_{\tau \overline{R\sig}} v^{- a_i \gamma_i b} \\
  &=& \sum_{b + \tau R\sig} \ti{\chi}((b + \tau R\sig)^{Q^f}) v^{a_i \gamma_i b} \1z_{\tau \overline{R\sig}} v^{- a_i \gamma_i b} \\
  &=& \sum_{b + \tau R\sig} 
  \ti{\chi}((a_i \gamma_i + \tau R\sig)^{Q^f})^{-1} X(a_i \gamma_i b + \tau R\sig) v^{a_i \gamma_i b} \1z_{\tau \overline{R\sig}} v^{- a_i \gamma_i b} \\
  &=& \ti{\chi}((a_i \gamma_i + \tau R\sig)^{Q^f})^{-1} x_X \overset{\eqref{ag^Q}}{=} x_X.
\eglnoz
As $\mu_{a_i \gamma_i}$ also fixes all the remaining summands of $w_\chi$ (see \eqref{w}), we conclude that $\mu_{a_i \gamma_i}(w_\chi) = w_\chi$. So for every generator $\gamma_i$, there exists an element $a_i$ in $\Fz_q\reg$ such that $\mu_{a_i \gamma_i}$ fixes $w_\chi$ for all $\chi$ in $\dualF \setminus \gekl{1}$. But the elements $\gekl{a_i \gamma_i}$ are again free generators of a subgroup $\Gamma'$ of $(K\sig)\reg$ with ${(K\sig)}^{\gekl{w}} = \Fz_q \times \Gamma'$.

The following lemma is proven with the help of the Pimsner-Voiculescu sequence, similarly to the proof of Lemma~\ref{mu_g=id}:
\blemma
\label{map-K->inj}
Assume that we can find unitaries $w_\chi$, $\chi \in \dualF \setminus \gekl{1}$, of the form \eqref{w} and free generators $\gamma_i$ for a choice of $\Gamma'$ (a subgroup of $(K\sig)\reg$ with $(K\sig)\reg = \Fz_q\reg \times \spkl{\tau} \times \Gamma'$) such that every $\mu_{\gamma_i}$ fixes $w_\chi$ for every $\chi \in \dualF \setminus \gekl{1}$. Set $\Gamma_i \defeq \gekl{\tau, \gamma_1, \dotsc, \gamma_{i-1}}$. Then the canonical map
\bgloz
  C_0(\Az_T\sig) \rtimes K\sig \rtimes (\Fz_q\reg \times \Gamma_i) \to C_0(\Az_T\sig) \rtimes K\sig \rtimes (\Fz_q\reg \times \Gamma_{i+1})
\egloz
induces an injective map on K-theory for every $i \in \Zz \nneg$.
\elemma
Now we know that
\bgloz
  C_0(\Az_T\sig) \rtimes K\sig \rtimes (K\sig)\reg \cong \ilim_i \gekl{C_0(\Az_T\sig) \rtimes K\sig \rtimes (\Fz_q\reg \times \Gamma_i)}.
\egloz
Thus the previous lemma, together with continuity of K-theory, implies that the canonical homomorphism
\bgloz
  C_0(\Az_T\sig) \rtimes K\sig \rtimes (\Fz_q\reg \times \spkl{\tau}) \to C_0(\Az_T\sig) \rtimes K\sig \rtimes (K\sig)\reg
\egloz
induces an injective map on K-theory. The reason is that, as it was explained above, we can choose $\Gamma'$ with free generators $\gamma_i$ for $\Gamma'$ and unitaries $w_\chi$ for $\chi \in \dualF \setminus \gekl{1}$ such that $\mu_{\gamma_i}(w_\chi)=w_\chi$ for all $\chi$ in $\dualF \setminus \gekl{1}$ as required by Lemma~\ref{map-K->inj}.

This observation, in combination with \eqref{AKK-Mor-AKsigKsig}, \eqref{A-Mor-infadele} and Proposition~\ref{K(AKFtau)} leads to the following
\bprop
\label{torsion-exists}
Assume that the inertia degree $f$ at infinity of $K / \Fz_q(T)$ satisfies $(f,q-1)=1$. Let $R$ be the ring of integers in $K / \Fz_p(T)$. Set $Q^f = \tfrac{q^f-1}{q-1}$. Then $K_*(\fA[R])$ contains $Q^f$-torsion.
\eprop

Now let us come to the proof of Theorem~\ref{mainthm2}. Let $K_1 / \Fz_p(T)$, $K_2 / \Fz_p(T)$ be function fields, both having the constant field $\Fz_q$. Assume that both of the function fields have only one single infinite place and let $f_i$ be the inertia degree at infinity of $K_i / \Fz_q(T)$. Suppose that $(f_i,q-1)=1$ for $i=1,2$. Let $R_1$ and $R_2$ denote the rings of integers in $K_1 / \Fz_p(T)$ and $K_2 / \Fz_p(T)$. In this situation, we can prove
\btheo
If $K_*(\fA[R_1]) \cong K_*(\fA[R_2])$, then the inertia degrees at infinity of $K_1 / \Fz_q(T)$ and $K_2 / \Fz_q(T)$ must coincide, i.e. $f_1 = f_2$.
\etheo
\bproof
Let $f_i$ denote the inertia degrees, for $i=1,2$. Since $K_*(\fA[R_1]) \cong K_*(\fA[R_2])$, we conclude by Proposition~\ref{torsion-exists} that $K_*(\fA[R_2])$ contains a torsion element of order $Q^{f_1}$. Hence, Corollary~\ref{torsion-order} tells us that every prime divisor of $Q^{f_1}$ divides $Q^{f_2}$. As the situation is symmetric in $K_1 / \Fz_q(T)$ and $K_2 / \Fz_q(T)$, we can equally well deduce that every prime divisor of $Q^{f_2}$ divides $Q^{f_1}$. Thus we conclude that $Q^{f_1}$ and $Q^{f_2}$ have the same set of prime divisors. The assertion of the theorem now follows from the following
\blemma
Let $q$ be an integer with $q>1$. Let $f_1$, $f_2$ be positive integers. Set $Q^{f_i} = \tfrac{q^{f_i}-1}{q-1}$ for $i=1,2$. If $Q^{f_1}$ and $Q^{f_2}$ have the same set of prime divisors, then $f_1=f_2$.
\elemma
\bproof
If $Q^{f_1}$ and $Q^{f_2}$ have the same prime divisors, then the set of their prime divisors is the set of prime divisors of their least common multiple $(Q^{f_1},Q^{f_2})$. We first prove
\bgl
\label{lcm}
  (Q^{f_1},Q^{f_2}) = Q^{(f_1,f_2)}.
\egl
Write $f=(f_1,f_2)$ and $f_1=fg_1$, $f_2=fg_2$ with $(g_1,g_2)=1$. Then
\bgloz
  Q^{f_1} = Q^f \cdot (1+q^f+\dotsb+(q^f)^{g_1-1})
\egloz
and
\bgloz
  Q^{f_2} = Q^f \cdot (1+q^f+\dotsb+(q^f)^{g_2-1}).
\egloz
We have to show $(1+q^f+\dotsb+(q^f)^{g_1-1},1+q^f+\dotsb+(q^f)^{g_2-1})=1$. As $(g_1,g_2)=1$, there are $a, b \in \Zz \nneg$ such that $ag_1=1+bg_2$. Thus
\bglnoz
  && (1+q^f+\dotsb+(q^f)^{g_1-1}) \cdot (1+(q^f)^{g_1}+\dotsb+((q^f)^{g_1})^{a-1}) \\
  &=& 1+q^f+\dotsb+(q^f)^{ag_1-1} = 1+q^f+\dotsb+(q^f)^{bg_2-1}+(q^f)^{bg_2} \\
  &=& (1+q^f+\dotsb+(q^f)^{g_2-1}) \cdot (1+(q^f)^{g_2}+\dotsb+((q^f)^{g_2})^{b-1}) + (q^f)^{bg_2}.
\eglnoz
Thus $(1+q^f+\dotsb+(q^f)^{g_1-1},1+q^f+\dotsb+(q^f)^{g_2-1})$ divides $(q^f)^{bg_2}$. If $b=0$, then $g_1$ has to be $1$ and we deduce $(1+q^f+\dotsb+(q^f)^{g_1-1},1+q^f+\dotsb+(q^f)^{g_2-1})=1$. If $b>0$, then $(1+q^f+\dotsb+(q^f)^{g_1-1},(q^f)^{bg_2})=1$ since every prime number $\pi$ that divides $(q^f)^{bg_2}$ also divides $q^f$, thus $\pi$ cannot divide $1+q^f+\dotsb+(q^f)^{g_1-1}$. This proves \eqref{lcm}.

Therefore, we can without loss of generality assume that $f_2$ divides $f_1$, i.e. we have $f_1 = f_2 f$. (In general, we can then apply our result to $f_1$ and $f=(f_1,f_2)$ and to $f_2$ and $(f_1,f_2)$.) Now consider the prime decompositions
\bgloz
  Q^{f_2} = \pi_1^{n_1} \dotsm \pi_r^{n_r} \text{ and } Q^{f_1} = \pi_1^{N_1} \dotsm \pi_r^{N_r}
\egloz
where $n_i$ and $N_i$ are positive integers. Here, we have used our assumption that the set of prime divisors of $Q^{f_1}$ and $Q^{f_2}$ coincide. Moreover, since $f_2$ divides $f_1$, we must have $N_i \geq n_i$ for all $1 \leq i \leq r$. It follows that
\bgloz
  q^{f_2}-1 = (q-1)\pi_1^{n_1} \dotsm \pi_r^{n_r} \text{ and thus } q^{f_2} = 1+(q-1)\pi_1^{n_1} \dotsm \pi_r^{n_r}.
\egloz
Analogously, we obtain $q^{f_1} = 1+(q-1)\pi_1^{N_1} \dotsm \pi_r^{N_r}$. But we also have
\bglnoz
  q^{f_1} &=& (q^{f_2})^f = (1+(q-1)\pi_1^{n_1} \dotsm \pi_r^{n_r})^f \\
  &=& \sum_{e=0}^f \binom{f}{e} (q-1)^e \pi_1^{n_1e} \dotsm \pi_r^{n_re} \\
  &=& 1 + (q-1)\sum_{e=1}^f \binom{f}{e} (q-1)^{e-1} \pi_1^{n_1e} \dotsm \pi_r^{n_re}.
\eglnoz
It follows that
\bgl
\label{comp-v0}
  \pi_1^{N_1} \dotsm \pi_r^{N_r} = \sum_{e=1}^f \binom{f}{e} (q-1)^{e-1} \pi_1^{n_1e} \dotsm \pi_r^{n_re}.
\egl
We now claim that for all $1 \leq i \leq r$, we have for all $2 \leq e \leq f$
\bgl
\label{v_pi}
  v_{\pi_i}(f \pi_1^{n_1} \dotsm \pi_r^{n_r}) < v_{\pi_i} (\binom{f}{e} (q-1)^{e-1} \pi_1^{n_1e} \dotsm \pi_r^{n_re}).
\egl
Here $v_{\pi_i}(x)$ for an integer $x$ is defined as $\max (\menge{\nu \in \Zz \nneg}{\pi_i^\nu \tei x})$.

To prove this, we compute $v_{\pi_i}(f \pi_1^{n_1} \dotsm \pi_r^{n_r}) = v_{\pi_i}(f) + n_i$ and
\bgln
\label{comp-v1}
  && v_{\pi_i}(\binom{f}{e} (q-1)^{e-1} \pi_1^{n_1e} \dotsm \pi_r^{n_re}) \\
  &=& v_{\pi_i}(f \dotsm (f-e+1)) - v_{\pi_i}(e!) + (e-1)v_{\pi_i}(q-1) + n_ie \nonumber \\
  &\geq& v_{\pi_i}(f) - v_{\pi_i}(e!) + (e-1)v_{\pi_i}(q-1) + n_ie \nonumber
\egln
Now expand $e$ with respect to $\pi_i$, i.e. write
\bgloz
  e = a_0 + a_1 \pi_i + \dotsb + a_m \pi_i^m \text{ with } a_0, \dotsc, a_m \in \gekl{0, \dotsc, \pi_i-1}.
\egloz
By \cite{Neu}, Chapter~II, Lemma~(5.6), we have
\bgl
\label{comp-v2}
  v_{\pi_i}(e!) = \tfrac{1}{\pi_i-1} \rukl{(\pi_i^m-1)a_m + (\pi_i^{m-1}-1)a_{m-1} + \dotsb + a_1}.
\egl
If $\pi_i = 2$, then $2 \tei Q^{f_1}$ implies that $q$ is odd, hence for every $2 \leq e$, we must have $(e-1) v_{\pi_i}(q-1) > 0$. Inserting this into \eqref{comp-v1}, we obtain
\bgl
\label{comp-v3}
  v_{\pi_i}(\binom{f}{e} (q-1)^{e-1} \pi_1^{n_1e} \dotsm \pi_r^{n_re}) > v_{\pi_i}(f) - v_{\pi_i}(e!) + n_ie.
\egl
Moreover, we always have by \eqref{comp-v2} that
\bgln
\label{comp-v4}
  v_{\pi_i}(e!) &\leq& (\pi_i^m-1)a_m + (\pi_i^{m-1}-1)a_{m-1} + \dotsb + a_1 \\
  &=& e - (a_m + a_{m-1} + \dotsb + a_1 + a_0) \leq e-1. \nonumber
\egln
Combining \eqref{comp-v3} and \eqref{comp-v4}, we obtain in the case $\pi_i = 2$:
\bgloz
  v_{\pi_i}(\binom{f}{e} (q-1)^{e-1} \pi_1^{n_1e} \dotsm \pi_r^{n_re}) > v_{\pi_i}(f) - (e-1) + n_ie.
\egloz
But this inequality also holds for $\pi_i > 2$, since we always have
\bgl
\label{comp-v5}
  v_{\pi_i}(\binom{f}{e} (q-1)^{e-1} \pi_1^{n_1e} \dotsm \pi_r^{n_re}) \geq v_{\pi_i}(f) - v_{\pi_i}(e!) + n_ie
\egl
by \eqref{comp-v1}. Moreover, if $\pi_i > 2$, then we deduce from \eqref{comp-v2} that
\bglnoz
  v_{\pi_i}(e!) &<& (\pi_i^m-1)a_m + (\pi_i^{m-1}-1)a_{m-1} + \dotsb + a_1 \\
  &=& e - (a_m + a_{m-1} + \dotsb + a_1 + a_0) \leq e-1.
\eglnoz
Plugging this inequality into \eqref{comp-v5}, we obtain
\bgloz
  v_{\pi_i}(\binom{f}{e} (q-1)^{e-1} \pi_1^{n_1e} \dotsm \pi_r^{n_re}) > v_{\pi_i}(f) - (e-1) + n_ie.
\egloz
Thus we have for every $\pi_i$:
\bglnoz
  && v_{\pi_i}(\binom{f}{e} (q-1)^{e-1} \pi_1^{n_1e} \dotsm \pi_r^{n_re}) > v_{\pi_i}(f) - (e-1) + n_ie \\
  &=& v_{\pi_i}(f) + n_i + n_i(e-1) - (e-1) = v_{\pi_i}(f) + n_i + (n_i-1)(e-1) \\
  &\geq& v_{\pi_i}(f) + n_i = v_{\pi_i}(f \pi_1^{n_1} \dotsm \pi_r^{n_r}).
\eglnoz
This proves \eqref{v_pi}. As a consequence, we deduce for all $1 \leq i \leq r$
\bgloz
  v_{\pi_i}(\sum_{e=1}^f \binom{f}{e} (q-1)^{e-1} \pi_1^{n_1e} \dotsm \pi_r^{n_re}) = v_{\pi_i}(f \pi_1^{n_1} \dotsm \pi_r^{n_r}) = v_{\pi_i}(f)+n_i.
\egloz
Comparing this with \eqref{comp-v0}, we conclude that
\bgloz
  N_i = v_{\pi_i}(\pi_1^{N_1} \dotsm \pi_r^{N_r}) = v_{\pi_i}(f)+n_i,
\egloz
thus $v_{\pi_i}(f)=N_i - n_i$ for all $1 \leq i \leq r$. This means that $\pi_1^{N_1-n_1} \dotsm \pi_r^{N_r-n_r}$ divides $f$. But $\pi_1^{N_1-n_1} \dotsm \pi_r^{N_r-n_r} = Q^{f_1} / Q^{f_2} = 1 + q^{f_2} + (q^{f_2})^2 + \dotsb + (q^{f_2})^{f-1} > f$ if $f>1$ because $q>1$. Thus we obtain a contradiction unless $f=1$, i.e. $f_1=f_2$.
\eproof
This completes the proof of Theorem~\ref{mainthm2}.
\eproof

Let us comment on the conditions in our main results. First of all, since our function field $K / \Fz_p(T)$ has only one single infinite place, the inertia degree $f$ at infinity of $K / \Fz_q(T)$ always divides the degree $n = [K : \Fz_q(T)]$. Thus the primeness assumption $(f,q-1)=1$ in Theorem~\ref{mainthm2} can be replaced by the stronger statement $(n,q-1)=1$. As the latter condition just requires to look at the global fields instead of the local fields at infinity, it might be simpler to verify.

Moreover, we can explicitly construct function fields $K / \Fz_p(T)$ satisfying all the conditions needed in Theorem~\ref{mainthm1}.

\bremark
\label{ex-ff}
For every prime power $q$ and every positive integer $f$, there exists a finite separable field extension $K / \Fz_p(T)$ such that
\begin{itemize}
\item
  the field of constants of $K / \Fz_p(T)$ is $\Fz_q$
\item
  $K / \Fz_p(T)$ has only one single infinite place
\item
  the inertia degree at infinity of $K / \Fz_q(T)$ is $f$.
\end{itemize}
\eremark

One possible construction goes as follows: Let $\zeta$ be a primitive ($q^f-1$)-th root of unity in $(\Fz_q)^{\text{alg}}$, so that $\Fz_{q^f} = \Fz_q(\zeta)$. Let $g(X) \in \Fz_q[X]$ be the minimal polynomial of $\zeta$ over $\Fz_q$. Set $K \defeq \Fz_q(X)$ and embed $\Fz_p(T)$ into $K$ via $T \ma g(X)^{-1}$. The corresponding extension $K / \Fz_p(T)$ is separable as $\Fz_q$ is perfect, so that the derivative of $g(X)$ does not vanish. Moreover, $K / \Fz_p(T)$ satisfies all the three conditions listed above by construction.

Remark~\ref{ex-ff}, in combination with Theorem~\ref{mainthm2}, shows that unlike in the case of number fields, there are infinitely many possibilities for the K-groups of ring C*-algebras attached to rings of integers in function fields. The reason why this is not the case for number fields is that there are not too many possibilities for the inertia degree at infinity in the number field case: A real place is said to have inertia degree $1$, and a complex place is said to have inertia degree $2$ (compare \cite{Neu}, Chapter~III, \S~1).

\end{document}